December 14, 2021

# New series identities with Cauchy, Stirling, and harmonic numbers, and Laguerre polynomials

Khristo N. Boyadzhiev
Ohio Northern University
Department of Mathematics and Statistics
Ada, Ohio 45810, USA
e-mail: k-boyadzhiev@onu.edu

**Abstract** In this article we use an interplay between Newton series and binomial formulas in order to generate a number of series identities involving Cauchy numbers (also known as Gregory coefficients or Bernoulli numbers of the second kind), harmonic numbers, Laguerre polynomials, and Stirling numbers of the first kind.

2010 Mathematics Subject Classification: Primary 11B65; Secondary, 05A19, 33B99, 30B10.

Keywords: Cauchy number, harmonic number, Stirling number, Newton series, binomial identity, Riemann's zeta function, Laguerre polynomial, series with special numbers.

## 1. Introduction

The Newton interpolation series is a classic tool in analysis with important applications ([12], [17], [18], [19]). In this paper, we use the Newton series in order to construct two summation formulas; the first one involving Cauchy numbers and the other one - Stirling numbers of the first kind. We obtain a number of series identities - closed form evaluations of several series involving Cauchy numbers together with harmonic numbers. We also present similar series with harmonic numbers and Stirling numbers of the first kind. The results are close to some classical results about such series.

The Cauchy numbers $c_n$ (see OEIS A006232 and A006233) are defined by the generating function

$$\frac{x}{\ln(x+1)} = \sum_{n=0}^{\infty} \frac{c_n}{n!} x^n \quad (|x|<1)$$

and they can also be defines as an integral of the falling factorial (see equation (2) below). These numbers are called Cauchy numbers of the first kind by Comtet [11]. The numbers $c_n / n!$ are called Bernoulli numbers of the second kind by some authors (see the comments in [3]). In Milne-Thomson's book [17] the Cauchy numbers appear as $B_n^{(n)}(1)$ on page 135. Blagouchine in his recent papers, [2], [3], and [4] used the notation $G_n = c_n / n!$ and the name Gregory coefficients for $G_n$.

A detailed study of the Cauchy numbers $c_n$ and several interesting identities involving these numbers were presented by Merlini et al in [16]. Cauchy numbers were also studied by Liu et al. [15] and Zhao [25].

The Cauchy numbers appear in some important formulas (like the Laplace summation formula in [17]) where they act like a counterpart to the Bernoulli numbers in the Euler-Maclaurin formula.

The Stirling numbers of the first kind $s(n,k)$ (OEIS A048994 and A008275) play a major role in this paper. They appeared for the first time in the works of James Stirling together with the Stirling numbers of the second kind $S(n,k)$ (for historical remarks see [6] and [22]). They are related to the Cauchy numbers by the formula ([11], [14])

$$c_n = \sum_{k=0}^{n} \frac{s(n,k)}{k+1} .$$

A detailed study of the Stirling numbers can be found in the books of Jordan [14], Comtet [11], and Graham et al. [13] (where the unsigned Stirling numbers of the first kind $(-1)^{n-k} s(n,k)$ are considered). Informative comments and interesting series containing $s(n,k)$ can be found in in Adamchik's paper [1] and in Section 2 of Blagouchine's paper [2]. Various results for series with Stirling numbers and harmonic numbers $H_n$ were recently obtained by Choi [10] and Wang and Lyu [23].

To give the reader an idea of our identities, here are three samples

$$\ln(1+x) = \sum_{n=0}^{\infty} \frac{(-1)^n c_n x^{n+1}}{(1+x)(1+2x)...(1+nx)}$$

$$2\zeta(3) - 1 = \sum_{n=1}^{\infty} \frac{(-1)^{n-1} c_n (H_n^2 + H_n^{(2)})}{n!n}$$

$$\sum_{j=1}^{k} \zeta(j+1) = \sum_{n=k}^{\infty} \frac{(-1)^{n-k} s(n,k) H_n}{(n+1)!} .$$

Starting from the Newton interpolation series with simple nodes $n = 0, 1, 2, ...$ we derive two special representations (Proposition 1 in Section 2 and Proposition 2 in Section 4). Using these two propositions and with the help of several binomial identities we generate a number of expansions involving Cauchy, harmonic, and Stirling numbers. Our results are presented in the form of examples. In Section 3 we work mostly with Cauchy numbers and in Section 4 we present series with Stirling numbers of the first kind and harmonic numbers. Section 5 contains a formula for series with Stirling numbers and Laguerre polynomials. Particular cases of this formula provide series with Laguerre polynomials and harmonic numbers.

## 2. The Newton interpolation series

We will use finite differences. Given an appropriate function $f(x)$, by definition we have

$$\Delta f(x) = f(x+1) - f(x),\ \Delta^n f(x) = \Delta(\Delta^{n-1} f(x)) \text{ for } n = 1, 2, \ldots \text{ and } \Delta^0 f(x) = f(x).$$

The Newton interpolation series for a given function $f(x)$ is the representation of this function in terms of the polynomials

$$(z)_n = z(z-1)\ldots(z-n+1),\ (z)_0 = 1$$

(see [6], [12], [17], [18]). Namely,

$$f(z) = \sum_{n=0}^{\infty} \frac{\Delta^n f(0)}{n!} z(z-1)\ldots(z-n+1). \tag{1}$$

This is the case of simple nodes $0,1,2\ldots$ . The Newton series expansion is not as popular as the Taylor series. One reason for this is the difficulty to study the remainder and to determine the area of convergence. There are other difficulties too. For instance, if we apply (1) to the function $f(x) = \sin(2\pi x)$, then $\Delta^n f(0) = 0$ for $n \ge 0$ and the series is identically zero.

Alexander Gelfond has done a thorough study of Newton interpolation series in chapter 2 of his book [12]. Another study can be found in chapter 5 of Norlund's book [18] (see also [19] and pp. 302-304 in Milne-Thomson's book [17]). The sufficient conditions for a function $f(z)$ to be represented by such a series are too elaborate to list here, however, we can say that if the function is analytic in a region of the form $\mathrm{Re}\,(z) > \lambda$ for some $\lambda$ and has moderate growth in that region, then the representation (1) holds for $\mathrm{Re}\,(z) > \lambda$. Moreover, for any two numbers $\varepsilon, R > 0$ the series is uniformly convergent for $\mathrm{Re}(z) \ge \lambda + \varepsilon, |z - \lambda| < R$. Boas and Buck show in [5], p. 34, that entire functions of exponential type less than $\ln 2$ have such representation.

The polynomials $(z)_n$ often appear in mathematics. For example, they are the ordinary generating functions for the Stirling numbers of the first kind $s(n,k)$ ([11] and [13])

$$z(z-1)\ldots(z-n+1) = \sum_{k=0}^{n} s(n,k)\, z^k .$$

They also appear in the definition of the Cauchy numbers $c_n$, $n = 0,1,\ldots$ ([11])

(2) $$c_n = \int_0^1 z(z-1)...(z-n+1)\,dz$$

with $c_0 = 1,\ c_1 = \frac{1}{2},\ c_2 = -\frac{1}{6},\ c_3 = \frac{1}{4}, ...$ etc. For more details on these numbers see Merlini et al. [16] and also Liu et al. [15].

Integrating the series (1) term by term and using equation (2) we come to the most interesting formula

(3) $$\int_0^1 f(x)\,dx = \sum_{n=0}^{\infty} \frac{c_n}{n!} \Delta^n f(0)$$

which clearly demonstrates the importance of the Cauchy numbers. This representation appears as equation (1) on page 277 in Jordan's book [14]. Note that in his book Jordan works with the numbers

$$b_n = \int_0^1 \binom{x}{n} dx = \frac{c_n}{n!}.$$

(A simpler version was used by Gregory in the 17[th] century.)

It is easy to compute that

(4) $$\Delta^n f(x) = (-1)^n \sum_{k=0}^{n} \binom{n}{k} (-1)^k f(x+k)$$

and in particular

$$\Delta^n f(0) = (-1)^n \sum_{k=0}^{n} \binom{n}{k} (-1)^k f(k)$$

Replacing the finite differences by the binomial expression in (3) we come to the following proposition.

**Proposition 1.** *For appropriate functions* $f(x)$ *we have the representation*

(5) $$\int_0^1 f(x)\,dx = \sum_{n=0}^{\infty} \frac{(-1)^n c_n}{n!} \left\{ \sum_{k=0}^{n} \binom{n}{k} (-1)^k f(k) \right\} .$$

As we shall see, this formula is a valuable tool for obtaining identities with Cauchy numbers.

**Definition**. The binomial transform of one sequence $a_0, a_1, ...$ is the new sequence defined by

$$b_n = \sum_{k=0}^{n} \binom{n}{k} (-1)^k a_k \quad \text{with inversion} \quad a_n = \sum_{k=0}^{n} \binom{n}{k} (-1)^k b_k .$$

The binomial transform of the sequence $\{f(k)\}$ is usually easier to compute compared to the computation of $\Delta^n f(0)$ so that formula (5) often works better than (3). A table of binomial transform formulas can be found in [6]. Using such formulas we will generate a number of examples in the next section.

Here is one simple demonstration how (5) works. Taking the function $f(z)=z^p$, where $p\geq 0$ is an integer we find by (5)

$$\frac{1}{p+1}=\sum_{n=0}^{\infty}\frac{(-1)^n c_n}{n!}\left\{\sum_{k=0}^{n}\binom{n}{k}(-1)^k k^p\right\}$$

and we recall that the Stirling numbers of the second kind $S(p,n)$ can be represented as

$$n!S(p,n)=(-1)^n\sum_{k=0}^{n}\binom{n}{k}(-1)^k k^p$$

so we get the identity

$$\frac{1}{p+1}=\sum_{n=0}^{p}c_n\, S(p,n)$$

(The infinite series terminates, as $S(p,n)=0$ for $n>p$). This is the Stirling inverse of the representation of the Cauchy numbers in terms of Stirling numbers of the first kind ([11], p, 294 and also [14], p. 267)

$$c_n=\sum_{p=0}^{n}\frac{s(n,p)}{p+1}.$$

## 3. Series with Cauchy numbers and harmonic numbers

In this section we use Proposition 1 together with various binomial identities to generate series identities involving Cauchy numbers and harmonic numbers.

**Example 1**. Let $y>0$ be arbitrary. We shall use the function

$$f(z)=\frac{y}{y+z},\ \ \operatorname{Re}(z)+y>0$$

together with the binomial identity (Entry (8.28) in [7])

$$\sum_{k=0}^{n}\binom{n}{k}(-1)^k\frac{y}{y+k}=\binom{n+y}{n}^{-1}.$$

Equation (5) implies

$$y\ln\left(1+\frac{1}{y}\right)=\sum_{n=0}^{\infty}\frac{(-1)^n c_n}{n!}\binom{n+y}{n}^{-1}=\sum_{n=0}^{\infty}\frac{(-1)^n c_n}{(y+1)(y+2)...(y+n)}$$

or

$$\ln\left(1+\frac{1}{y}\right)=\sum_{n=0}^{\infty}\frac{(-1)^n c_n}{y(y+1)(y+2)...(y+n)}$$

(cf. [24] and for a similar representation see also p. 244 in [19]).

Replacing $1/y$ by $x$ we find the remarkable expansion of $\ln(1+x)$

$$(6)\qquad \ln(1+x)=\sum_{n=0}^{\infty}\frac{(-1)^n c_n x^{n+1}}{(1+x)(1+2x)...(1+nx)}\ .$$

When $x=1$ this gives

$$\ln 2=\sum_{n=0}^{\infty}\frac{(-1)^n c_n}{(n+1)!}.$$

Using the asymptotic behavior of the Cauchy numbers at infinity ([11], p. 294)

$$\frac{c_n}{n!}\sim\frac{(-1)^{n+1}}{n\ln^2 n}\quad (n\to\infty)$$

it is easy to verify (by the ratio test, for instance) that the series in (6) is convergent for all $x>0$.

Before continuing further we recall that the harmonic numbers $H_n$ are defined for $n=0,1,...$ by

$$H_n=1+\frac{1}{2}+...+\frac{1}{n}\quad (n>0),\ H_0=0$$

and it is often convenient to express them in the form

$$H_n=\psi(n+1)+\gamma$$

where $\psi(z)=\frac{d}{dz}\ln\Gamma(z)$ is the digamma function and $\gamma=-\psi(1)$ is Euler's constant. This formula provides an extension of the harmonic numbers as a holomorphic function $H_z$ on $\operatorname{Re}(z)>-1$. We use also the generalized harmonic numbers $H_n^{(p)}$, where $H_0^{(p)}=0$ and for $n>0$

$$H_n^{(p)}=1+\frac{1}{2^p}+...+\frac{1}{n^p}=\zeta(p)+\frac{(-1)^{p-1}}{(p-1)!}\psi^{(p-1)}(n+1)\quad (p>1).$$

**Example 2.** Let $m > 0$ be an integer. We take the function $f(z) = \dfrac{1}{(z+m)^2}$ and use the binomial identity

$$\sum_{k=0}^{n} \binom{n}{k} (-1)^k \frac{1}{(k+m)^2} = \frac{(m-1)!n!(H_{n+m} - H_{m-1})}{(n+m)!}$$

(Entry (8.38) in [7]). Equation (5) implies

$$\text{(7)} \qquad \frac{1}{(m+1)!} = \sum_{n=0}^{\infty} \frac{(-1)^n c_n (H_{n+m} - H_{m-1})}{(n+m)!}.$$

For $m = 1, m = 2,$ and $m = 3$ we find correspondingly

$$\frac{1}{2} = \sum_{n=0}^{\infty} \frac{(-1)^n c_n H_{n+1}}{(n+1)!}$$

$$\frac{1}{6} = \sum_{n=0}^{\infty} \frac{(-1)^n c_n (H_{n+2} - 1)}{(n+2)!}$$

$$\frac{1}{12} = \sum_{n=0}^{\infty} \frac{(-1)^n c_n (2H_{n+3} - 3)}{(n+3)!}$$

etc.

**Example 3**. Here we exploit the binomial formula (Entry (9.16) in [7])

$$\text{(8)} \qquad \sum_{k=0}^{n} \binom{n}{k} (-1)^k H_k^{(3)} = -\frac{1}{2n} (H_n^2 + H_n^{(2)})$$

($n = 1, 2...$ ) where

$$H_k^{(3)} = 1 + \frac{1}{2^3} + ... \frac{1}{k^3} = \zeta(3) + \frac{1}{2} \psi''(k+1).$$

We apply (5) to the function $f(z) = \zeta(3) + \frac{1}{2} \psi''(z+1)$. In this case

$$\int_0^1 \zeta(3) + \frac{1}{2} \psi''(z+1)\, dz = \zeta(3) + \frac{1}{2} \psi'(z+1)\Big|_0^1 = \zeta(3) + \frac{1}{2} (\psi'(2) - \psi'(1)) = \zeta(3) - \frac{1}{2}$$

and therefore,

$$2\zeta(3) - 1 = \sum_{n=0}^{\infty} \frac{(-1)^n c_n}{n!} \left\{ 2 \sum_{k=0}^{n} \binom{n}{k} (-1)^k H_k^{(3)} \right\}$$

$$= \sum_{n=1}^{\infty} \frac{(-1)^n c_n}{n!} \left\{ \sum_{k=0}^{n} \binom{n}{k} (-1)^k H_k^{(3)} \right\} = \sum_{n=1}^{\infty} \frac{(-1)^{n-1} c_n (H_n^2 + H_n^{(2)})}{n!n}$$

(the term for $n = 0$ in the first sum is zero and we can start the summation from $n = 1$). Thus

$$(9) \qquad 2\zeta(3) - 1 = \sum_{n=1}^{\infty} \frac{(-1)^{n-1} c_n (H_n^2 + H_n^{(2)})}{n!n}.$$

It can be seen from the digamma definition of the harmonic numbers that

$$\lim_{p \to 0} \frac{1}{p} (H_p^2 + H_p^{(2)}) = 2\zeta(3)$$

which can be considered the zero term in the sum. With this agreement, starting the summation from $n = 0$ we can write

$$-1 = \sum_{n=0}^{\infty} \frac{(-1)^{n-1} c_n (H_n^2 + H_n^{(2)})}{n!n}.$$

**Example 4**. For completeness we will reprove here two known identities. First, the classical representation

$$(10) \qquad \gamma = \sum_{n=0}^{\infty} \frac{(-1)^{n-1} c_n}{n!n}$$

obtained by Lorenzo Mascheroni in 1790. Informative historical notes on this series can be found on page 406 in Blagouchine's paper [2].

We take here $f(z) = \psi(z+1) + \gamma$, $f(k) = H_k$ and use the simple identity (9.3a) in [7]

$$\sum_{k=0}^{n} \binom{n}{k} (-1)^k H_k = -\frac{1}{n}.$$

Clearly

$$\int_0^1 \psi(z+1) + \gamma dz = \ln \Gamma(z+1) \Big|_0^1 + \gamma = \gamma$$

and (5) implies (9) immediately.

Next we prove the representation

$$(11) \qquad \frac{\pi^2}{6} - 1 = \sum_{n=1}^{\infty} \frac{(-1)^{n-1} c_n H_n}{n!n}$$

(see [16] and also equation (44) on p. 413 in [2]). We take $f(z) = \zeta(2) - \psi'(z+1)$ so that $f(k) = H_k^{(2)} = 1 + \frac{1}{2^2} + ... + \frac{1}{k^2}$, $f(0) = 0$ . We need also the identity ((9.4b) in [7])

$$\sum_{k=0}^{n} \binom{n}{k} (-1)^k H_k^{(2)} = -\frac{H_n}{n}.$$

Then we compute

$$\int_0^1 \zeta(2) - \psi'(z+1)\, dz = \zeta(2) - \psi(z+1)\Big|_0^1 = \zeta(2) - \psi(2) + \psi(1) = \zeta(2) - 1 .$$

and (11) follows from (5).

Like in the previous example we notice that

$$(12) \qquad \lim_{p \to 0} \frac{H_p}{p} = \lim_{p \to 0} \frac{\psi(p+1) + \gamma}{p} = \zeta(2) = \frac{\pi^2}{6}$$

and if we agree to write

$$\frac{\pi^2}{6} = \frac{H_n}{n}\bigg|_{n=0}$$

equation (11) becomes (with summation from $n = 0$)

$$-1 = \sum_{n=0}^{\infty} \frac{(-1)^{n-1} c_n H_n}{n!n}.$$

Merlini, Sprugnoli and Verri [16] found (10) and (11) by using a symbolic Laplace summation formula. Blagouchine derived (11) from a special series transformation formula (p. 412 in [2]). With the help of his formula he obtained also several more interesting representations of this kind. Two finite sums connecting Cauchy numbers and harmonic numbers can be found in [25].

**Example 5.** In this example we use the constant

$$M_1 = \int_1^2 \frac{\psi(1+t) + \gamma}{t}\, dt \approx 0.86062$$

and we will prove that

$$(13) \qquad M_1 = \frac{1}{2} + \sum_{n=0}^{\infty} \frac{(-1)^{n-1} c_n H_n}{(n+1)!}.$$

For this purpose we take the function $f(z) = \frac{\psi(z+1) + \gamma}{z+1}$ and apply the binomial formula

(14) $$\sum_{k=0}^{n}\binom{n}{k}(-1)^k \frac{H_k}{k+1} = -\frac{H_n}{n+1}$$

(entry (9.32) from [7]). The left hand side in (5) becomes

$$\int_0^1 \frac{\psi(z+1)+\gamma}{z+1}dz = \int_1^2 \frac{\psi(x)+\gamma}{x}dx = \int_1^2 \frac{\psi(x+1)-1/x+\gamma}{x}dx$$

$$= M_1 - \int_1^2 \frac{1}{x^2}dx = M_1 - \frac{1}{2}$$

and therefore,

$$M_1 - \frac{1}{2} = \sum_{n=0}^{\infty} \frac{(-1)^{n-1} c_n H_n}{(n+1)!}.$$

The constant $M_1$ is important, because it represents the value of some interesting series. It was shown in [8] that

$$M_1 = \sum_{n=1}^{\infty} \frac{1}{n} \ln\left(1 + \frac{1}{n+1}\right)$$

$$M_1 = \sum_{n=1}^{\infty} \frac{(-1)^{n-1}}{n}\left(n - \zeta(2) - \zeta(3) - \ldots - \zeta(n)\right)$$

(*the first term in this series is* 1). Also

$$M_1 = \sum_{n=1}^{\infty} H_n^{-}\left(\zeta(n+1) - 1\right)$$

where

$$H_n^{-} = 1 - \frac{1}{2} + \frac{1}{3} + \ldots + \frac{(-1)^{n-1}}{n}$$

are the skew harmonic numbers. Now we have one more series in this list.

## 4. Series with Stirling numbers of the first kind and harmonic numbers

In this section we assume that $f(z)$ is a function analytic in some half plane $\mathrm{Re}(z) > \lambda$, where $\lambda < 0$, and has moderate growth. Being analytic in a neighborhood of the origin, $f(z)$ will be represented by a Taylor series centered at $z = 0$.

Now we will manipulate the representation (1) in a different way. We replace $(z)_n$ by

$$z(z-1)...(z-n+1)=\sum_{m=0}^{n} s(n,m)\, z^m$$

and change the order of summation to get

$$f(z)=\sum_{m=0}^{\infty} z^m \left\{\sum_{n=0}^{\infty} \frac{s(n,m)}{n!}\Delta^n f(0)\right\}.$$

Comparing this to the Taylor series

$$f(z)=\sum_{m=0}^{\infty} z^m \frac{f^{(m)}(0)}{m!}$$

we conclude that

$$\frac{f^{(m)}(0)}{m!}=\sum_{n=0}^{\infty}\frac{s(n,m)}{n!}\Delta^n f(0)\,.$$

Applying equation (4) now we come to an important identity.

**Proposition 2**. *Under the initial assumption on the function $f(z)$, for every $m \geq 0$ we have the representation*

$$\text{(15)} \qquad \frac{f^{(m)}(0)}{m!}=\sum_{n=0}^{\infty}\frac{(-1)^n s(n,m)}{n!}\left\{\sum_{k=0}^{n}\binom{n}{k}(-1)^k f(k)\right\}.$$

(The summation, in fact, starts from $n=m$ since $s(n,m)=0$ for $n<m$ .)

When $m=1$ we have $s(n,1)=(-1)^{n-1}(n-1)!$ and the formula takes the form

$$f'(0)=\sum_{n=1}^{\infty}\frac{1}{n}\left\{\sum_{k=0}^{n}\binom{n}{k}(-1)^{k-1} f(k)\right\}.$$

This evaluation was obtained by a different method in [9].

Formula (15) is an efficient mechanism for producing various series with Stirling numbers of the first kind together with special numbers. It resembles formula (5), where the Cauchy numbers are replaced by Stirling numbers and on the left hand side we have now derivatives. Most examples from the previous section can be repeated with the same function $f(z)$ and the same binomial identity, but now using formula (15) instead of (5).

**Example 6**. First we prove a classical result. Taking $f(z)=\psi(z+1)+\gamma$ , $f(k)=H_k$ we have

$$\sum_{k=0}^{n}\binom{n}{k}(-1)^k f(k)=-\frac{1}{n}$$

(see Example 4). The Taylor series of $\psi(z+1)+\gamma$ is well-known

(16) $$\psi(z+1)+\gamma=\sum_{k=1}^{\infty}(-1)^{k-1}\zeta(k+1)\,z^k, \quad |z|<1$$

and (15) produces the representation ( $k\geq 1$ )

$$\zeta(k+1)=\sum_{n=k}^{\infty}\frac{(-1)^{n-k}\,s(n,k)}{n!n}.$$

This is almost a folklore mathematical result. In a more general form it can be found on p. 166 in Jordan's book [14]. Comments and extensions can be found in Adamchik's paper [1], Section 5, and also in Blagouchine's paper [2], on p. 412.

Note that the numbers

$$\begin{bmatrix} n \\ k \end{bmatrix}=(-1)^{n-k}\,s(n,k)$$

are called the unsigned Stirling numbers of the first kind. They are thoroughly studied in the classical book of Graham, Knuth, and Patashnik [13]. These numbers will appear in the examples that follow.

**Example 7**. Working with the function $f(z)=\zeta(2)-\psi\,'(z+1)$ , $f(k)=H_k^{(2)}$ (as in Example 4) we have for $n\geq 1$

$$\sum_{k=0}^{n}\binom{n}{k}(-1)^k\,f(k)=-\frac{H_n}{n}$$

and also from (16)

$$\zeta(2)-\psi\,'(z+1)=\zeta(2)+\sum_{k=0}^{\infty}(-1)^{k+1}(k+1)\zeta(k+2)\,z^k$$

$$=\sum_{k=1}^{\infty}(-1)^{k+1}(k+1)\zeta(k+2)\,z^k\,.$$

This way for $k\geq 1$ we have

(17) $$\zeta(k+2)=\frac{1}{k+1}\sum_{n=k}^{\infty}\frac{(-1)^{n-k}\,s(n,k)H_n}{n!n}\,.$$

The equality holds also for $k=0$ under the agreement (12).

Continuing further in this direction we can use identity (8) from Example 3 with the same function $f(z)=\zeta(3)+\frac{1}{2}\psi''(z+1), f(k)=H_k^{(3)}$. Again from (16)

$$f(z)=\frac{1}{2}\sum_{k=1}^{\infty}(-1)^{k-1}(k+1)(k+2)\zeta(k+3)z^k$$

and (15) implies

$$\text{(18)}\qquad \zeta(k+3)=\frac{1}{(k+1)(k+2)}\sum_{n=k}^{\infty}\frac{(-1)^{n-k}s(n,k)}{n!n}\left(H_n^2+H_n^{(2)}\right)$$

for any $k\geq 1$.

The identities (17) and (18) were found independently by a different method by Wang and Lyu in [23] (see p.171 there).

**Example 8.** Let again $k\geq 1$. Applying (15) to the function

$$f(z)=\frac{\psi(z+1)+\gamma}{z+1}, \text{Re}(z)>-1$$

and the binomial identity (14) we find

$$\text{(19)}\qquad \sum_{j=1}^{k}\zeta(j+1)=\sum_{n=k}^{\infty}\frac{(-1)^{n-k}s(n,k)H_n}{(n+1)!}$$

(see [21, p. 28]). To prove this we use the Taylor series (16) and notice that for $|z|<1$

$$\frac{\psi(z+1)+\gamma}{z+1}=\frac{1}{z+1}\sum_{k=1}^{\infty}(-1)^{k-1}\zeta(k+1)z^k=\sum_{k=1}^{\infty}\left\{(-1)^{k-1}\sum_{j=1}^{k}\zeta(j+1)\right\}z^k$$

and (15) implies (19) immediately.

**Example 9.** Here we use the identity

$$\sum_{k=0}^{n}\binom{n}{k}(-1)^k\frac{(-1)^k}{(k+1)^2}=\frac{H_{n+1}}{n+1}$$

together with the function $f(z)=(z+1)^{-2}, \text{Re}(z)>-1$. For $|z|<1$ we have the Taylor series

$$\frac{1}{(z+1)^2}=\sum_{k=0}^{\infty}(-1)^k(k+1)z^k.$$

Therefore, from (15) we find the curious companion to equation (19)

$$\text{(20)}\qquad k+1=\sum_{n=k}^{\infty}\frac{(-1)^{n-k}s(n,k)H_{n+1}}{(n+1)!}$$

([21, p. 29]). More generally, using the identity from Example 2

$$\sum_{k=0}^{n} \binom{n}{k} (-1)^k \frac{1}{(k+m)^2} = \frac{(m-1)!\,n!\,(H_{n+m} - H_{m-1})}{(n+m)!}$$

with the function

$$f(z) = (z+m)^{-2} = \sum_{k=0}^{\infty} \frac{(-1)^k (k+1) z^k}{m^{k+2}}$$

we find

$$\frac{k+1}{m^{k+2}(m-1)!} = \sum_{n=k}^{\infty} \frac{(-1)^{n-k} s(n,k) (H_{n+m} - H_{m-1})}{(n+m)!}. \tag{21}$$

where $m = 1$ gives equation (20). In the same line, using the identity

$$\sum_{k=0}^{n} \binom{n}{k} (-1)^k \frac{1}{(k+m)^3} = \frac{(m-1)!\,n!}{2(n+m)!} \left[ (H_{n+m} - H_{m-1})^2 + H_{n+m}^{(2)} - H_{m-1}^{(2)} \right]$$

((8.42) in [7]) together with the function

$$f(z) = (z+m)^{-3} = \sum_{k=1}^{\infty} \frac{(-1)^k (k+1)(k+2) z^k}{m^{k+3}}$$

we come to the series

$$\frac{(k+1)(k+2)}{m^{k+3}(m-1)!} = \sum_{n=k}^{\infty} \frac{(-1)^{n-k} s(n,k)}{2(n+m)!} \left[ (H_{n+m} - H_{m-1})^2 + H_{n+m}^{(2)} - H_{m-1}^{(2)} \right]. \tag{22}$$

**Example 10**. Here we use the identity

$$\sum_{k=0}^{n} \binom{n}{k} (-1)^k H_{k+1} = -\frac{1}{n(n+1)}$$

together with the function $f(z) = \psi(z+2) + \gamma$, $f(k) = H_{k+1}$. We write

$$f(z) = \psi(z+2) + \gamma = \psi(z+1+1) + \gamma = \frac{1}{1+z} + \psi(z+1) + \gamma$$

for $\operatorname{Re}(z) > -1$. From (16)

$$f(z) = \sum_{k=1}^{\infty} (-1)^{k-1} (\zeta(k+1) - 1)\, z^k, \quad |z| < 1$$

and (15) produces the representation ([14, p. 339])

(23) $$\zeta(k+1)-1=\sum_{n=k}^{\infty}\frac{(-1)^{n-k}s(n,k)}{(n+1)!\,n} \quad (k\ge 1)\,.$$

**Example 11**. Now we present a series connecting the Stirling numbers $s(n,k)$ and the central binomial coefficients $\binom{2n}{n}$. Starting from the identity (Entry (8.43) in [7]).

$$\sum_{k=0}^{n}\binom{n}{k}(-1)^k\frac{1}{2k+1}=\frac{4^n}{2n+1}\binom{2n}{n}^{-1}$$

we introduce the function $f(z)=(2z+1)^{-1}$, $\mathrm{Re}(z)>-\frac{1}{2}$ with Taylor series for $|z|<1/2$

$$\frac{1}{2z+1}=\sum_{k=0}^{\infty}(-1)^k 2^k z^k\,.$$

From (15) we have the strange representation

(24) $$2^k=\sum_{n=k}^{\infty}\frac{(-1)^{n-k}4^n s(n,k)}{n!(2n+1)}\binom{2n}{n}^{-1}\,.$$

All these series with positive terms are very slowly convergent.

## 5. Series with Laguerre polynomials and harmonic numbers

In this section we present series with the classical Laguerre polynomials $L_n(x)$ which have the representation

$$L_n(x)=\sum_{k=0}^{n}\binom{n}{k}(-1)^k\frac{x^k}{k!}\,.$$

Let $x>0$ and consider the function $f(z)=\dfrac{x^z}{\Gamma(z+1)}$ so that $f(k)=\dfrac{x^k}{k!}$. Proposition 2 implies the formula

(25) $$\left(\frac{d}{dz}\right)^m\frac{x^z}{\Gamma(z+1)}\bigg|_{z=0}=\sum_{n=m}^{\infty}\frac{(-1)^n s(n,m)}{n!}L_n(x)\,.$$

This way for every $m=1,\,2,\ldots$, we have different series with Laguerre polynomials.

**Example 12.** When $m=1$, $(-1)^n s(n,1) = -(n-1)!$ and we come to the known series

$$-\gamma - \ln x = \sum_{n=1}^{\infty} \frac{1}{n} L_n(x) . \tag{26}$$

This is entry 5.11.1 (3) in the reference book [20]. Interestingly, this equation resembles equation (101) in [4].

For $m=2$ we have $(-1)^n s(n,2) = (n-1)! H_{n-1}$ and this provides the new series

$$\frac{\gamma^2}{2} + \gamma \ln x - \frac{\pi^2}{12} + \frac{\ln^2 x}{2} = \sum_{n=1}^{\infty} \frac{1}{n} H_{n-1} L_n(x) . \tag{27}$$

For $m>2$ the derivatives in (25) become very long. For $m=3$ we just write

$$\left(\frac{d}{dz}\right)^3 \frac{x^z}{\Gamma(z+1)} \Bigg|_{z=0} = \frac{1}{2} \sum_{n=k}^{\infty} \frac{1}{n} \left( H_{n-1}^2 - H_{n-1}^{(2)} \right) L_n(x) \tag{28}$$

etc.

## References

[1] **Victor Adamchik,** On Stirling numbers and Euler sums, *J. Comput. Applied Math.*, 79 (1997), 119-130.

[2] **Iaroslav V. Blagouchine**, Two series expansions for the logarithm of the gamma function involving Stirling numbers and containing only rational coefficients for certain arguments related to $\pi^{-1}$, *J. Math. Anal. Appl.*, 442 (2016), 404–434.

[3] **Iaroslav V. Blagouchine,** A note on some recent results for the Bernoulli numbers of the second kind, *J. Integer Seq*. 20 (2017), Article 17.3.8.

[4] **Iaroslav V. Blagouchine**, Three notes on Ser's and Hasse's representatations for the zeta-function, *INTEGERS*, 18A (2018).

[5] **Ralph P. Boas, JR., and R. Creighton Buck**, *Polynomial Expansions of Analytic Functions*, Academic Press, 1964.

[6] **Khristo N. Boyadzhiev**, Close encounters with the Stirling numbers of the second kind, *Math. Mag*., 85, No. 4, (2012), 252-266.

[7] **Khristo N. Boyadzhiev**, *Notes on the Binomial Transform, Theory and Table*, World Scientific, 2018.

[8] **Khristo N. Boyadzhiev,** A special constant and series with zeta values and harmonic numbers, *Gazeta Matematica, Seria A*, vol. 115 (3-4) (2018), 1-16.

[9] **Khristo N. Boyadzhiev**, Evaluation of series with binomial sums, *Anal. Math*., 40 (1), (2014), 13-23

[10] **Junesang Choi,** Certain identities involving Stirling numbers of the first kind, *Far East J. Math. Sci*., 103(2) (2018), 523-539.

[11] **Louis Comtet**, *Advanced Combinatorics*, Kluer, 1974.

[12] **Alexander O. Gelfond**, *Calculus of finite differences*. Translated from the Russian. International Monographs on Advanced Mathematics and Physics. Hindustan Publishing Corp., Delhi, 1971. (Russian original *Ischislenie konechnih raznostei,* Moskva, 1959.)

[13] **Ronald L. Graham, Donald E. Knuth, Oren Patashnik**, *Concrete Mathematics*, Addison-Wesley Publ. Co., New York, 1994.

[14] **Charles Jordan**, *Calculus of finite differences*, Chelsea, New York, 1950 (First edition: Budapest 1939)

[15] **Hong-Mei Liu, Shu-Hua Qi, and Shu-Yan Ding**, Some recurrence relations for Cauchy numbers of the first kind, *J. Integer Seq*., 13 (2010), Article 10.3.8

[16] **Donatella Merlini, Renzo Spignoli, M. Cecilia Verri**, The Cauchy numbers, *Discrete Math*., 306 (2006), 1906-1920.

[17] **Louis M. Milne-Thomson**, *The Calculus of Finite Differences*, MacMillan, 1951.

[18] **Niels E. Nørlund,** *Leçons Sur les Séries d'Interpolation,* Gauthier-Villars, 1926.

[19] **Niels E. Nørlund,** *Vorlesungen Uber Differenzenrechnung*, Chelsea, 1954.

[20] **Anatoli P. Prudnikov, Yuri A. Brychkov, Oleg I. Marichev,** I*ntegrals and Series, Vol.2, Special Functions*, Gordon and Breach 1986.

[21] **Jefferey O. Shallit**, Metric theory of Pierce expansions, *Fibonacci Quart*., 24, (1986), 22-40.

[22] **Ian Tweddle,** *James Stirling's Methodus Differentialis: An Annotated Translation of Stirling's Text*, Springer, 2003.

[23] **Weiping Wang and Yanhong Lyu**, Euler sums and Stirling sums, *J. Number Theory*, 185 (2018), 160-193.

[24] **George N. Watson**, The transformation of an asymptotic series into a convergent series of inverse factorials, *Rend. Circ. Mat. Palermo* 34 (1912), 41-88.

[25] **Feng-Zhen Zhao**, Sums of products of Cauchy numbers, *Discrete Math*., 309 (2009), 3830-3042.